\documentclass[12pt]{article}
\usepackage{color}
\usepackage{amsfonts,amssymb,epsfig,amscd,amsmath,marvosym,graphics}
\usepackage[latin1]{inputenc}
\usepackage{setspace}
\usepackage[english]{babel}
\newtheorem{thm}{Theorem}[section]
\newtheorem{cor}{Corollary}[section]
\newtheorem{lem}{Lemma}[section]
\newtheorem{pro}{Proposition}[section]
\newtheorem{dfn}{Definition}[section]
\newtheorem{rmq}{Remark}[section]
\newtheorem{expl}{Example}[section]

\input{amssym.def}
\input{amssym.tex}
\definecolor{couleurliensref}{rgb}{0.,0.,1.}
\definecolor{couleurliens}{rgb}{0.,0.,0.}
\definecolor{couleurliensurl}{rgb}{.0,.0,.0}
 \usepackage[pdftex,colorlinks=true,urlcolor=couleurliensurl,linkcolor=couleurliens,citecolor=couleurliensref]{hyperref}
 \def\dessous#1\sous#2{\mathrel{\mathop{\kern0pt#2}\limits_{#1}}}

\newcommand{\R}{\mathbb R} 
 
\newcommand{\f}{\mathcal{F}}

\newcommand{\beq}{\begin{eqnarray}}
\newcommand{\eeq}{\end{eqnarray}}
\newcommand{\bpro}{\begin{pro}}
\newcommand{\epro}{\end{pro}}
\newcommand{\blem}{\begin{lem}}
\newcommand{\elem}{\end{lem}}
\newcommand{\bdfn}{\begin{dfn}}
\newcommand{\edfn}{\end{dfn}}
\newcommand{\bcor}{\begin{cor}}
\newcommand{\ecor}{\end{cor}}
\newcommand{\bthm}{\begin{thm}}
\newcommand{\ethm}{\end{thm}}
\newcommand{\bex}{\begin{expl}}
\newcommand{\eex}{\end{expl}}
\newcommand{\brmq}{\begin{rmq}}
\newcommand{\ermq}{\end{rmq}}
\newcommand{\benum}{\begin{enumerate}}
\newcommand{\eenum}{\end{enumerate}}
\newcommand{\bitem}{\begin{itemize}}
\newcommand{\eitem}{\end{itemize}}

\title{ Hamiltonian mean curvature flow
\footnotetext{This work has been completed during the first author' stay at ICTP within STEP program.}\\[.2cm] }
 \author{ 
Djid\'em\`e F. Hou\'enou\footnote{Institut de Math\'ematiques et de Sciences-Physiques (IMSP)} and 
  L\'eonard Todjihound\'e$^*$\\[.2cm] 
}
\date {}

\begin{document}
\maketitle
\begin{abstract}
 
 \noindent Let $(\Sigma,\omega)$ be a compact Riemann surface with constant curvature $c$. In this work, 
we proved that the mean curvature flow of a given Hamiltonian diffeomorphism on $\Sigma$ 
provides a smooth path in  $Ham(\Sigma)$, the group of all Hamiltonian diffeomorphisms of $\Sigma$.
 This result gives a proof, in the case of graph of Hamiltonian diffeomorphisms to the conjecture 
of Thomas and Yau asserting that the mean curvature flow of a compact embedded Lagrangian 
submanifold $S$ with zero Maslov class in a Calabi-Yau manifolds $M$ exists for all time and converges 
smoothly to a special Lagrangian submanifold in the Hamiltonian isotopy class of $S$.
\end{abstract}




\section{ Introduction } 

The deformation of maps between Riemannian manifolds has been intensively studied for a long time. 
The idea is to find a natural process to deform a map to a \textsl{canonical} one. The harmonic 
heat flow is probably the famous example although Ricci flow and mean curvature flow are also 
well used. The latter is an evolution process under which a submanifold of a given manifold evolves in the 
direction of its mean curvature vector. From the first variation formula for the volume functional, 
one easily observes that the mean curvature flow represents the most effective way to decrease the 
volume of a submanifold such that it is very useful  when minimal submanifolds or volume minimizer
 submanifolds are sorted for under suitable conditions. Many results have been found for mean curvature 
flow in codimension one while the higher codimension is still receiving attention of number of researchers. 

The simplest case of higher codimension is the mean curvature flow of surface in 4-dimensional manifolds. 
 It compounds two important classes known as \textsl{symplectic mean curvature flow} and 
\textsl{Lagrangian mean curvature flow} ; since being symplectic or Lagrangian is preserved along the flow.
 In the last decades, several works in geometric analysis field research are devoted to these classes 
(see e.g \cite{xh-jl, ks, mk1, mk2}).
 
It is well known that the geometric structures of the ambient space plays a fundamental  role  when  
studying the  existence and the properties of the mean curvature flow. For instance, to the authors knowledge, 
symplectic mean curvature flow exists only when the ambient manifold carries at least an (almost) K\"ahler-Einstein 
structure 
\cite{mk2} ; or particularly  when it is Calabi-Yau (target spaces for superstring compactification). 
Moreover Lagrangian and special Lagrangian submanifolds of Calabi-Yau manifolds are considered as the 
cornerstones for understanding the mirror symmetry phenomenon between pairs of Calabi-Yau manifold both of 
the categorical point of view and from a physical-geometrical standpoint (see e.g \cite{ssz}). 

There exists a  cohomology class attached to any given Lagrangian submanifold of a symplectic manifold : the  
{\sl Maslov class}. This class can be represented by a closed 1-form expressed solely 
in term of the mean curvature of the submanifold and the symplectic form of the ambient manifold. 
Therefore one observes that minimal Lagrangian  has zero Maslov class. In light of this fact  we consider 
the deformation of Hamiltonian diffeomorphism by the mean curvature flow. We will call this flow, when it exists in 
the group of Hamiltonian diffeomorphisms, a {\it Hamiltonian mean curvature flow }.  

Let us recall that being graph and Lagrangian are preserved along the mean curvature flow (see e.g. \cite{mk2},
 \cite{mtw1} ) and these two properties together yield the mean curvature flow of symplectomorphism
under the hypothesis that  the universal covering of the ambient manifold is of the type 
$\mathbb S^2\times \mathbb S^2$, $\R^2\times\R^2$ or $\mathbb H^2\times\mathbb H^2$.

Later on, in \cite[Theorem 1.1]{mtw3}, the author proved the following : 

\noindent{\bf Theorem}$ $

{\it Let $(\Sigma_1, \omega_1)$  and $(\Sigma_2, \omega_2)$ be two homeomorphic compact Riemann surface of the 
same constant curvature $c =-1$, 0, or 1. Suppose $\Sigma$ is the graph of a symplectomorphism 
$f : \Sigma_1 \longrightarrow \Sigma_2$ as a Lagrangian submanifold of 
$M = (\Sigma_1\times \Sigma_2, \omega_1-\omega_2)$ and $\Sigma_t$ is the mean curvature flow with initial
surface $\Sigma_0=\Sigma$. Then $\Sigma_t$ remains the graph of a symplectomorphism $f_t$ along the mean 
curvature flow. The flow exists smoothly for all time and $\Sigma_t$ converges smoothly to a minimal 
Lagrangian submanifold as $t\longrightarrow \infty$.}

In account of this we proved that  deforming Hamiltonian diffeomorphism by the mean curvature flow provides 
a path in the group of  Hamiltonian diffeomorphisms. The result is stated as follows :
 
\bthm $ $

Let $(\Sigma,\omega)$ be a compact connected Riemann surface with non-negative constant curvature.
Then any Hamiltonian diffeomorphism on $\Sigma$ deforms through Hamiltonian diffeomorphisms under 
the mean curvature flow.
\ethm

In the sequel we set some of the needed materials in the approach of Hamiltonian mean curvature flow 
and recall the technical tools used for the study of such geometric evolution problem. We then
 discuss the Hamiltonian property of the time slice of the flow.
\section{ Preliminaries }
Throughout this exposition, all manifolds are smooth and  closed (compact without boundary) unless it is 
stated otherwise. 

\bdfn $ $

Let $M$ be a differentiable manifold.
The mean curvature motion of $S\subset M$ is a 1-parameter family of immersions of  submanifolds $S_t$ 
in $M$ which admits a parametrization  $F_t:S\longmapsto S_t\subset M$ over $S$ with normal velocity equal
to the mean curvature vector i.e. 
\beq\label{e1}
\left(\frac{\partial}{\partial t}F_t(x)\right)^\bot&=& H(F_t(x)), \qquad x\in S ,\\\nonumber
 F_0&=&id.
\eeq
\edfn

\brmq $ $

The mean curvature motion is a non linear weakly parabolic system for $F$ and is invariant under 
reparametrization of $S$. Indeed, by coupling with a diffeomorphism $\varphi$ of $S$, the flow can be made 
into a normal direction, i.e. 
\beq\label{e2}
 \frac{\partial}{\partial t}F_t\big(\varphi(x)\big)= H\big(F_t(\varphi(x))\big).
\eeq
For any smooth compact initial data, one can establish the short time existence solution for (\ref{e2}) 
and the uniqueness of the solution for suitable conditions on the initial data.
\ermq

Let $(\Sigma,\omega)$ be a compact Riemann surface with a constant curvature $c$ and 
let $f\in\mathcal D{\rm iff}(\Sigma,\omega)$ be a  diffeomorphism of $\Sigma$. Put $M=\Sigma\times\Sigma$ 
and denote by $S$ the graph of $f$. 

The mean curvature flow of $f$ is realized through the mean curvature motion of $S$. In fact, knowing 
that the mean curvature flow preserved the graph and Lagrangian properties, one obtains the flow in the group 
of symplectic diffeomorphisms $S{\rm ymp}(\Sigma,\omega)$. 

\bdfn\label{d1} $ $

A diffeomorphism $f$ on a symplectic manifold $(\Sigma,\omega)$ is said to be Hamiltonian if there 
exists a smooth function $G: \Sigma\longrightarrow \R$ such that   $f\in \{f_s\}_s$,
 where $\{f_s\}_s$ is the Hamiltonian flow of $X=X_{G}$, i.e the family of  
diffeomorphisms obtained by solving the ordinary differential equation : 
\beq
  \left\{
\begin{array}{rcl}
 \frac{\partial}{\partial s}f_s (x) & =  & X (f_s(x)) \\[.2cm]\nonumber
 f_0(x) & = & x.
\end{array}
\right.
\eeq
\edfn

The definition above is a classical definition of Hamiltonian diffeomorphism. In this definition, the 
vector field does not depend on time and is often refer to as autonomous vector field. The analog
for time-dependent vector field is  the characterization  using the flux homomorphism. It is stated as
 in the following : 

\bdfn \label{d2}$ $

 The time-one map of a symplectic isotopy (path to identity) with null flux is called a  Hamiltonian 
diffeomorphism (see  \cite{ab1}).
\edfn
 
Therefore, the flux is an obstruction for symplectomorphism which is isotopic to identity to be a Hamiltonian 
diffeomorphism. Let us recall the definition of the flux homomorphism. For more details we refer 
to \cite{ab2}  where a comprehensive exposition is made.
The flux homomorphism is defined as follows, [Theorem 3.1.1, \cite{ab2}] :
\beq
\begin{array}{cccl}
\widetilde {Flux} : & \widetilde{S{\rm ymp}_0}(\Sigma,\omega)& \longrightarrow & H^1(\Sigma,\R) \\[.2cm]
      & \{\tilde f_s\} &\longmapsto & \left[\int_0^1 f_s^*(i_{X_s}\omega) \  ds \right]
\end{array}
\eeq
where $\widetilde{S{\rm ymp}_0}(\Sigma,\omega)$ is the universal covering of the connected component of the identity 
in the group of symplectomorphisms, $\{\tilde f_s\}$ is a homotopy class of an isotopy $f_s$ generated by
 $X_s$, and $f_s^*(i_{X_s}\omega)$ stands for the pull-back of the form 
$i_{X_s}\omega$ (interior product of $\omega$ by $X_s$). Notice that the flux homomorphism descends to the  
group ${S{\rm ymp}_0}(\Sigma,\omega)$ ; for convenience, we will recall this definition  in the following section.

Let $\{{f_s}\}_{0\leq s\leq 1}$ be a symplectic isotopy to $f$, generated by the time-dependent vector 
field $X_s$ and denote $\f$, its flux form i.e.
 $$\f =\int_0^1\!\! f_s^*i_{X_s}\omega\  ds.$$

The cohomology class $[\f]$ depends only on the homotopy classes of the isotopy $\{f_s\}_{0\leq s\leq 1}$
 relatively to fixed ends. The behavior of $\f$ under the mean curvature flow  will be the major ingredient 
for the purpose of preserving Hamiltonian condition by the mean curvature flow since it is well known that  
%
the mean curvature flow of symplectomorphism exists smoothly for all time and converges  
(see e.g. \cite{ks, mk2, mtw1, mtw2}).

In the sequel, we compute the evolution equation of $\f$ and use it to find out under 
which hypothesis the Hamiltonian property is preserved along the flow.

\section{Hamiltonian property along the flow}

Recall that the group of Hamiltonian diffeomorphisms $Ham(\Sigma,\omega)$ is the kernel of an onto homomorphism 
$$Flux :S{\rm ymp}_0 (\Sigma,\omega) \longrightarrow H^1(\Sigma, \R)/\Gamma_\omega$$
where the flux group $\Gamma_\omega=Flux\bigg(\pi_1 \Big(S{\rm ymp}_0(\Sigma,\omega)\Big)\bigg)$ is finitely 
generated but is not known to be discrete in all cases. 
Hence the most one can say in general is that $Ham(\Sigma,\omega)$ sits inside the identity component 
$S{\rm ymp}_0(\Sigma,\omega)$ as the leaf of a foliation. Therefore, we do not use the topology on 
$Ham(\Sigma,\omega)$ induced 
from $S{\rm ymp}_0(\Sigma,\omega)$ but instead use the topology on $Ham(\Sigma,\omega)$ induced from the 
$C^2$-topology on the Lie algebra of Hamiltonian functions with zero mean. Thus a neighborhood of the 
identity consists of all time 1-maps of Hamiltonian flows generated by Hamiltonians $G_s$ that are 
sufficiently small in the $C^2$-topology.

 Let $f:\Sigma\longrightarrow\Sigma$ be a Hamiltonian diffeomorphism on a compact Riemann surface with 
constant curvature. As asserted above, we have two different ways to regard $f$.  Let us consider a
symplectic isotopy view point, meaning  $f$ is the end point of some symplectic isotopy 
$\{f_s\}_{0\leq s\leq 1}$ with  zero flux. The mean curvature flow of $f$ gives rise to a 2-parameter family 
of symplectomorphisms  $\{f_{s,t}\}$ satisfying :
\beq (A)\ 
\left\{
\begin{array}{ccll}
f_{1,0} & = & f&   \\[.2cm]\nonumber 
f_{0,t}  =  id_{\Sigma}, & & f_{1,t}=f_t  \qquad \text{ for each t }&\\\\\nonumber 
Flux\{f_{s,0}\} & =& \left[\int_0^1 f_{s,0}^*(i_{X_{s,0}}\omega) ds\right] = [\f_0]= 0 \\\\\nonumber 
\frac{\partial}{\partial s} f_{s,t} & = & X_{s,t}  \\[.2cm]\nonumber
\frac{\partial}{\partial t} f_{s,t} & = & H_{s,t} \\[.2cm]\nonumber
\frac{\partial}{\partial t}X_{s,t} & = &\frac{\partial}{\partial s}H_{s,t} -[H_{s,t},X_{s,t}]\\[.2cm]
\end{array}\right.
\eeq
where $X_{s,t}$ is the isotopy vector field and  $H_{s,t}$ the mean curvature vector field ; 
$s$ is the isotopy parameter and $t$ is the one of the mean curvature flow.

\bdfn \label{d3}$ $

Let $X$ be a vector field with a local 1-parameter group $(\varphi_t)_t$ of local diffeomorphisms and 
$S$ a tensor field on a differentiable manifold $M$. The Lie derivative of $S$ in the direction
$X$ is defined as
\beq 
\varphi_t^*L_X S := \left(\frac{d}{dt}\varphi_t^* S\right)\Bigg|_{t=0} 
\eeq
\edfn

Let $M$ and $N$ be two differentiable manifolds. Assume $\varphi_t : M\longrightarrow N$ is a smooth 
1-parameter family of maps between $M$ and $N$ and $\omega_t$  is a smooth family of forms on $N$. 
Then $\varphi_t^*\omega_t$  is a smooth family of forms on $M$ and the basic formula of differential 
calculus of forms gives (see e.g \cite{mj}) :
\beq
\frac{d}{dt}\varphi_t^*\omega_t= \varphi_t^*L_{X_t}\omega_t + \varphi_t^*\frac{d}{dt}\omega_t,
\eeq
where $X_t$ is the tangent vector field along $\varphi_t$. 
\blem $ $

The flux form $\f_t$ of the isotopy $f_{s,t}$, satisfies the following equation :
\beq\label{e3}
 \frac{\partial}{\partial t} \f_t= f_t^*i_{H_t}\omega +d K_t
\eeq
where $K_t=\int_0^1f^*_{s,t} \omega(X_{s,t}, H_{s,t})ds$ and $f_t$ (the time $t$-slice of the flow) 
is the time-one map of the isotopy $\{f_{s,t}\}_{0\leq s\leq 1}$. 
\elem

{\it Proof} : 

Taking into account the fact that $X_{s,t}$ and $H_{s,t}$ are symplectic vector fields,
and using Definition \ref{d1}, a direct computation yields 
\beq
 \frac{\partial}{\partial t} \f_t & = &\int_0^1 \frac{\partial}{\partial t}\Big(f_{s,t}^*i_{X_{s,t}}\omega\Big)ds\\\nonumber
                                &=& \int_0^1\left( f_{s,t}^*L_{H_{s,t}}i_{X_{s,t}}\omega + 
                                   f_{s,t}^*\frac{\partial}{\partial t}i_{X_{s,t}}\omega \right)ds\\\nonumber
               &=&  \int_0^1\left( f_{s,t}^*di_{H_{s,t}}i_{X_{s,t}}\omega + 
                                   f_{s,t}^*\frac{\partial}{\partial s}i_{H_{s,t}}\omega -f_{s,t}^*i_{[H_{s,t},X_{s,t}]}\omega \right)ds\\\nonumber
&=&  \int_0^1\left( df_{s,t}^*i_{H_{s,t}}i_{X_{s,t}}\omega + 
                                   f_{s,t}^*\frac{\partial}{\partial s}i_{H_{s,t}}\omega +f_{s,t}^*L_{X_{s,t}} i_{H_{s,t}}\omega \right)ds\\\nonumber
&=&  \int_0^1\left( df_{s,t}^*i_{H_{s,t}}i_{X_{s,t}}\omega + 
                                  \frac{\partial}{\partial s}\left( f_{s,t}^*i_{X_{s,t}}\omega\right)  \right)ds\\\nonumber
                            & = &  f_t^*i_{H_t}\omega +d\int_0^1 f_{s,t}^*\omega(X_{s,t}, H_{s,t})ds
\eeq
$\hfill{\square}$

It was discovered in 1965 by V. P. Maslov that there is a cohomology class which appears naturally in 
the resolution by the Hamilton-Jacobi method of the Schr\"odinger equation of quantum physics. In mathematics, 
it is a cohomology class attached to a given Lagrangian submanifold a symplectic manifold.
This class is an important cohomology invariant and is called the Maslov class. Since J.-M. Morvan's work, 
it has been found possible to express this class solely in terms of the Riemmanian structure of the 
Lagrangian immersion associated to the K\"ahler metric on a symplectic manifold (see e.g \cite{dr}).

\bdfn $ $

Let $(M,\omega)$ be a K\"ahler-Einstein  2$n$-dimensional manifold and $L\hookrightarrow M$ be an  
immersed Lagrangian submanifold of $M$. Then the Maslov class of $L$ is defined by :
\beq
\frac{n}{\pi}[i_H\omega],
\eeq 
where $H$ is the mean curvature vector along $L$.
\edfn

Therefore we obtain the following : 
\blem \label{l2}$ $

Let $(f_s)_{0\leq s\leq 1}$ be a symplectic isotopy to a Hamiltonian diffeomorphism $f$. 
Then the cohomology class of the flux form of $f_s$ deforms to the Maslov class of $S$ 
by the mean curvature flow.
\elem

{\it Proof } : 

From equation (\ref{e3}) one gets :
 \beq
 \frac{\partial}{\partial t}[\f_t]=\left[\frac{\partial}{\partial t}\f_t\right]= [f_t^*i_{H_t}\omega].
\eeq$\hfill{\square}$

As an immediate consequence to Lemma \ref{l2}, the following holds :

\bpro\label{p2} $ $

Let $\Sigma$ be a compact connected Riemann surface with constant curvature and $f\in Ham(\Sigma)$. 
Suppose $S=graph f$  has zero  Maslov class, then the flux of any symplectic isotopy to $f$ 
is preserved along the mean curvature flow.
\epro

{\it Proof} : 

Let  $\{f_s\}_{0\leq s \leq 1}$ be a symplectic isotopy to $f$ generated by the vector field $X_s$. The mean 
curvature flow of $f$ is a 2-parameter family  $\{f_{s,t}\}_{s,t}$ of symplectomorphisms. So using 
Lemma \ref{l2} and taking into account the fact that $\Sigma$ is connected,
one concludes that the cohomology class of the flux form is constant along the flow.$\hfill{\square}$

We now state the main results of this work.

\bthm\label{t2} $ $

Let $\Sigma$ be a compact connected Riemann surface with constant curvature and $f\in Ham(\Sigma)$.
Assume that $S=graph f$ has zero Maslov class.
Then any Hamiltonian diffeomorphism on $\Sigma$ deforms through Hamiltonian diffeomorphisms by the mean
curvature flow.
\ethm
{\it Proof} : 

Let $f\in Ham(M,\omega)$. There exists a symplectic isotopy $\{f_s\}$ to $f$ such that 
the mean curvature flow of $f$ is a 2-parameter as in system $(A)$. We know that the flow exists 
(see e.g. \cite{mtw1}).
 So using the  Proposition \ref{p2}, for each time $t$, the flux of the isotopy $\{f_{s,t}\}$ to $f_t$ is zero.
Therefore  $f_t$ is a Hamiltonian diffeomorphism.$\hfill{\square}$

Thus we call {\sl Hamiltonain mean curvature flow} a mean curvature flow for which any time slice of the flow 
is Hamiltonian or equivalently a mean curvature flow of Lagrangian graphs Hamiltonian isotopic to the diagonal.

%

%

\subsection{Non-negative curvature case}
 In this section we assume that $\Sigma$ has  non-negative constant curvature and observe that the 
assumption of zero Maslov class can be removed. We have the following : 
\bthm \label{t3}$ $

Let $(\Sigma,\omega)$ be a compact connected Riemann surface with constant curvature $c$ and 
$f\in Ham(\Sigma)$. If $c$ is non-negative, then $f$ deforms through Hamiltonian diffeomorphisms under 
the mean curvature flow.
\ethm

{\it Proof} :

$M=\Sigma\times\Sigma$ is compact and its universal covering is either  
$\mathbb{S}^2\times \mathbb{S}^2$ or $\R^2\times \R^2$. The submanifold $S$ (graph of $f$) is Lagrangian 
w.r.t $\omega'=\omega\ominus \omega$ and symplectic w.r.t $\omega\oplus \omega$. Then $f$ deforms through 
symplectomorphism \cite{mtw1}. What is left is  to prove that each slice $f_t$ of the flow is Hamiltonian, 
i.e a time one map of some symplectic isotopy $\{f_{s,t}\}_{0\leq s\leq 1}$ with zero flux.
\benum
\item  Suppose  $\Sigma$ is  elliptic  ($c>0$), then $H^1(\Sigma,\R)$ is trivial and so is $H^1(S,\R)$ ;
 thus the flux is preserved. Since its initial value is zero, then each $t$-slice of the flow is Hamiltonian.

\item If $c=0$, then $M$ is Calabi-Yau. Let $\theta_t$ be the Lagrangian angle of $S_t$. 
The mean curvature form satisfies 
$$ i_{H_t}\omega = d\theta_t $$
which implies that the Maslov class vanishes.
Thus, the flux is preserved along the flow and since its initial value is zero, we deduce that each $f_t$ 
is Hamiltonian. $\hfill{\square}$
\eenum

\bdfn $ $

A diffeomorphism on $\Sigma$ is called a minimal diffeomorphism if its graph is a minimal embedding 
in $\Sigma\times\Sigma$.
\edfn

As a consequence of Theorem \ref{t3}, we obtain the following corollaries :

\bcor \label{c1}$ $

Let $(\Sigma,\omega)$ be a compact connected Riemann surface with non-negative constant curvature $c$ 
and $f\in Ham(\Sigma)$. As $t \longrightarrow \infty$, a sequence of the mean curvature flow of the 
graph of $f$ converges to a smooth minimal Hamiltonian graph.
\ecor

\bcor $ $ 

Let $(\Sigma,\omega)$ be a compact connected Riemann surface with strictly positive constant curvature $c$ 
and $f\in Ham(\Sigma)$. Then the Hamiltonian mean curvature fow of $S$=graph$f$ exists for all time $t$,
 each $S_t$ can be written as a graph of a Hamitonian diffeomorphism $f_t$.
 The sequence of submanifolds $S_t$ converges to the diagonal as $t$ goes to infinity.
\ecor

The proof of these corollaries are the same as in \cite{mtw1} ; in addition with the preserving Hamiltonian 
property from Theorem \ref{t3}. 
%
%
%
%
%

 A particular example of calibrated submanifolds was first introduced by Harvey and Lawson. These submanifolds 
are known as {\textsl special Lagrangian} (see definition bellow). It is not hard to check that calibrated submanifolds are volume 
minizer in their homology class  so special Lagrangian are minimal submanifolds.

\bdfn $ $

A Lagrangian submanifold in a Calabi-Yau manifold $(M,\Omega)$ is called special if it has constant 
Lagrangian angle.
\edfn


The Theorem \ref{t3} and its Corollary \ref{c1} give the proof in the case of graph of Hamiltonian 
diffeomorphism to the conjecture of Thomas and Yau  asserting that the mean curvature flows of 
a Lagrangian submanifold $S$ with zero Maslov class exits for all time and converges to a special Lagrangian
submanifold in the Hamiltonian isotopy class of $S$. We proved that  the Hamiltonian isotopy is nothing else 
but the path obtained  by the mean curvature flow.
 
\bthm $ $

Let $S$ be a graph of some Hamiltonian diffeomorphism $f$ on a flat torus $T^2$. Then the Hamiltonian 
mean curvature flow of $S$ exists for all time and converges to a special Hamiltonian graph isotopic to $S$.
\ethm

{\it Proof} :  

 $Ham(T^2)$ is contractible, so every $f\in Ham(T^2)$ flows through Hamiltonian 
diffeormorphisms to the identity  which graph (the diagonal in $T^2\times T^2$) is a minimal surface.
Then the mean curvature form is exact which infers that the Lagrangian angle is constant. Thus the 
limit is a special Lagrangian submanifold.$\hfill\square$

\section*{Acknowledgment :} 

 We would like to thank Professor A. Banyaga for pointing out this problem, 
 Professors J. Li and M-T. Wang for their very helpful suggestions and comments.
We are also grateful to Professor Claudio Arrezo for his invaluable suggestions.


\addcontentsline{toc}{chapter}{\textit{Bibliography}}{ }

\end{document}